
\def\Pn#1{{\bf P}^#1}
\def\slutt{\sqcap\!\!\!\!\sqcup}
\font\first=cmbx10 at 15pt

  \magnification=1200\parindent=0pt
\centerline {\first Conic bundles in projective fourspace}
\bigskip\bigskip
\centerline {Robert Braun and Kristian Ranestad}\bigskip\bigskip\bigskip

P. Ellia and G.Sacchiero have shown that if $S$ is a smooth
surface in $\Pn 4$ which is ruled in conics, then $S$ has degree 4 or 5
(cf.
[ES]).  In this paper we give a proof of this result combining the
ideas of Ellia and Sacchiero as they are used in the paper of the
second author on plane curve fibrations [Ra] and the recent work of G.
Fl\o ystad and the first author bounding the degree of smooth surfaces
in $\Pn 4$ not of general type [BF].\par

Let $S$ be a smooth conic bundle in $\Pn 4$. Let $V$ denote the
hypersurface
which is the union of the planes of the conics on $S$.
 Let $G\subset \Pn 9$ be the Grassmannian of planes in $\Pn 4$ in
the Pl{\"u}cker embedding.  Since the hypersurface $V$ contains a one
dimensional family of planes, we may associate a curve $C_V\subset G$ whose
points correspond to the planes in $V$.    In the natural incidence variety
in
$G\times \Pn 4$ between points and planes, there is a $\Pn 2$-fibration $W$
over
$C_V$  whose projection into $\Pn 4$ is $V$. If $C_V$ is not smooth
consider its
normalization $\widetilde{C_V}$, and the pullback $\tilde W$ of $W$ to
$\widetilde{C_V}$.   The strict transform $\tilde S$ of $S$ in $\tilde W$
is
clearly smooth, since the map $\tilde S\to S$ is birational. So on the
complement of some possible $(-1)$-curves this map is an isomorphism.  But
if
$C_V$ is singular there are two plane curves, possibly infinitely close,
which
are mapped into the same plane in $\Pn 4$, this is a contradiction.
Therefore
$C_V$ is smooth. Let $g$ be the genus of $C_V$ and let  $\delta$=deg$C_V$,
thus
$\delta$ is also the degree of the hypersurface $V$.

The proof is an exploitation of the relations between the invariants of $S$
and
$C_V$. On the one hand they combine with the results of [BF] to give the
upper
bound $d\leq 42$ for the degree of $S$.  On the other hand the curve $C_V$
in the
Grassmannian inside $\Pn 9$, has a genus which is high compared to the
degree.
Thus Castelnuovo bounds show that the span of this curve is a subspace of
$\Pn
9$.  We analyse the intersection of the linear span of $C_V$ and the
Grassmannian, which is a variety $T$ cut out by quadrics. The lines and
conics
on this variety give rise to special curves on $S$ which together with
the bound for the degree allow us to conclude.
 We work over an algebraically closed field of characteristic 0.

\proclaim (1.1) Lemma.  If $S$ has degree $d$ and sectional genus $\pi$,
 then $$3d\geq
4\delta,\leqno {(a)}$$ $$\pi -1=d+2g-2-\delta,\leqno {(b)}$$
$$d^2-9d-8(2g-2)+2\delta =0.\leqno {(c)}$$
$$\pi -1={d^2\over 8}-{d\over 8}-{3\delta\over 4}.\leqno {(d)}$$\par
$Proof$.
Part $(a)$ is the fact that there is a positive number of singularities for
the map $\pi: S\to C_V$:  If a fibre is nonreduced i.e. a double line, then
the
line would have nonpositive selfintersection on $S$, but every fibre is a
conic
in a plane so computing the arithmetic genus of a fibre in two ways we get
a
contradiction.   Therefore the number
$$c_2(\Omega_S-\pi^{\ast}\Omega_{C_V})= 3d-4\delta$$ counts a nonnegative
finite
number of singular points.\par Part $(b)$ is a straightforward calculation
using
adjunction on $W$. Part $(c)$ is the double point formula for smooth
surfaces in
$\Pn 4$ applied to $S$ (cf.[HR,p.434]). Part $(d)$ follows from $(b)$ and
$(c)$.
$\slutt$\medskip
 \proclaim (1.2) Remark.  From Severis theorem
it follows that the projection of $W$ into $\Pn 4$ is linearly normal.
 Thus by Riemann Roch $\delta \leq 2+3g$.\par
\proclaim (1.3) Remark.  Smooth surfaces on quadrics are well understood
and
those on cubics are classified recently  (cf. [A],[K]) so we only need to
worry
about $\delta\geq 4$.
  In fact there are conic bundles on quadrics, and any conic bundle on a
cubic
is also on a quadric.\par

\proclaim (1.4) Lemma.  $\delta \leq 3$ or
$$g-1\geq{1\over 9}\delta^2-{5\over 8}\delta.$$\par
$Proof.$  Since ${4\over 3} \delta -{9\over 2}\geq 0$ when $\delta \geq 4$,
 the relations (1.1c) and (1.1a) yields $$\eqalign
{0=&d^2-9d-8(2g-2)+2\delta\cr\geq& {16\over
9}\delta^2-12\delta-8(2g-2)+2\delta\cr} $$ from
which the lemma follows.$\slutt$\medskip This inequality beats the genus
bound
for curves in $\Pn 6$ (cf. [HJ]):
\proclaim (1.5) Proposition.
$C_V$ is contained in a $\Pn 5$, and if $C_V$ spans a $\Pn 5$,
then it lies on a surface of degree 4.\par $Proof.$
  If $C_V$ spans a $\Pn
6$, then the genus bound says that $${1\over 10}({\delta
}^2-7\delta+12)\geq
{1\over 9}\delta^2-{5\over 8}\delta+1,$$
 which yields $$\delta^2+{27\over 4}\delta - 18\leq
0.$$ i.e. $\delta\leq 2$.\par
If $C_V$ is not contained in a surface of degree 4 in $\Pn 5$,
then the genus formula (cf.[HJ]) yields $${1\over 10}(\delta^2-5\delta
+10)\geq
{1\over 9}\delta^2-{5\over 8}\delta +1.$$
Thus $${\delta }^2 -{45\over 4} \delta \leq 0.$$ Thus the proposition
follows from
\proclaim (1.6) Lemma.   If $C_V$ spans $\Pn 5$ and does not lie on a
surface
of degree 4, then $\delta \geq 12$.\par
$Proof$.  If $C_V$ is rational or elliptic, then its degree is at least 5
or 6,
while the above Castelnuovo bound
 says that $g\leq 7$ when $\delta \leq 11$. It remains only to check that
(1.1c) has no integral solutions,
 which is straightforward.$\slutt$
\bigskip

{\bf  2  Bound for the degree of $S$}\medskip
{}From the results of [BF] we can show
\proclaim (2.1) Proposition. If $S$ is not
on a cubic hypersurface, then $d\leq 42$.\par $Proof.$  We distinguish
between
the cases whether $S$ lies on a quartic or a quintic hypersurface or not,
and
apply Roths theorem [Ro] to study the genus of a general hyperplane
section.\par   Case 1:  Assume $S$ is not contained in a quintic
and $d>25$.
 Then by Roths theorem a general hyperplane section of $S$ is also not
contained
in a quintic (in $\Pn 3$).  Hence the genus bound for space curves (cf.
[GP])
gives $$\pi - 1\leq {d^2\over 12} +d.$$
Combined with (1.1d) and using $\delta \leq {3d\over 4}$ from (1.1a) we
find
$d\leq 40$ in this case.\medskip
Case 2:  Assume $S$ is contained in a quintic, not contained in a quartic
and
$d>17$.  As in [BF 1.1b] let
$$\gamma ={d^2\over 10}+{d\over 2}+1-{2r(5-r)\over 5}-\pi$$
where $0\leq r\leq 4$ and $d+r\equiv 0\quad ({\rm mod}\quad 5)$.  By the
genus
bound for space curves (cf. [GP]) $\gamma$ is a non-negative integer
 satisfying $$\pi -1\leq {d^2\over 10}+{d\over 2}-\gamma$$ and furthermore
(cf.[BF 1.1e] )$${d^3\over 150}-{d\over 6}\leq\chi ({\cal O}_S)+{\gamma
^2\over
2}+\gamma({d\over 5}+{5\over 2}).$$
The first inequality combined with (1.1d) leads to
$$\gamma \leq {d\over 80}(95-2d).$$
Hence a priori $d\leq 47$.  Moreover the maximal value of $\gamma$ in the
range $18\leq d\leq 47$ is 14.
Now (1.1c) combined with (1.1a) yields
$$ \chi ({\cal O}_S) =1-g=-{d^2\over 16}+{9d\over 16} -{\delta\over 8
}\leq -{d^2\over 16}+{9d\over 16}.$$
Inserting this and $\gamma=14$ in the second inequality gives
$${d^3\over 150}-{d\over 6}\leq -{d^2\over 16}+{9d\over 16}+{14d\over
5}+133.$$
Evaluating shows that $d\leq 30$ in this case.\medskip
Case 3:  Assume $S$ is contained in a quartic and $d>10$.  As in [BF 1.1b]
let
$$\gamma={d^2\over 8}+1-{3r(4-r)\over 8}-\pi$$
where $0\leq r\leq 3$ and $d+r\equiv 0 \quad ({\rm mod}\quad 4)$.   By the
genus
bound for space curves (cf. [GP])
$\gamma$ is a non-negative integer and $$\pi -1\leq {d^2\over 8}-\gamma.$$
Combined with (1.1d) and (1.1a) this gives
$$\gamma \leq {d\over 8}+{3\delta\over 4}\leq {11d\over 16}.$$
We have (cf. [BF 1.1e])
$${d^3\over 96}-{d^2\over 16}-{d\over 24}+{5\over 4}\leq\chi ({\cal
O}_S)+{\gamma
^2\over 2}+\gamma({d\over 4}+{3\over 2}).$$
Putting things together (as in case 2) leads to
$${d^3\over 96}-{209d^2\over 512}-{157d\over 96}+{5\over 4}\leq 0$$ which
yields
$d\leq 42$.$\slutt$
 \proclaim
(2.2) Corollary. If $S$ is not on a cubic, then $\delta\leq 31$.\par
$Proof$
Combine (2.1) with (1.1a).$\slutt$
\bigskip
{\bf  3  Some geometry of $C_V$}\medskip
Now $C_V$ is a curve on the Grassmannian G, which is a variety cut out by
quadrics.    If $L$ is the linear span of $C_V$, and $T$ is the
irreducible component of $G\cap L$ which contains $C_V$, then $T$ is a
quadric
or lies on more than one quadric in $L$.  In each possible case we may
describe
the family of planes in $\Pn 4$ corresponding to the  closed points of $T$.
\proclaim (3.1) Lemma.  The closed points on a line in $G$ correspond to
the
pencil of planes through a line in a $\Pn 3$.\smallskip The closed points
on a
conic in $G$ whose plane is not contained in $G$, correspond to the planes
of
one of the pencils of a quadric of rank 4 in $\Pn 4$.\par
$Proof.$  Easy.$\slutt$

\proclaim (3.2) Lemma. The curve $C_V$ is rational or elliptic,
or $T$ is a plane,  or a quadric surface in a $\Pn 3$ or the whole $\Pn 3$,

or a cubic scroll or a cone over a quartic curve or a delPezzo
in  $\Pn 4$, or a quadric hypersurface in $\Pn 4$, or $C_V$ spans a $\Pn
5$. \par
 $Proof.$ The linear span $L$ of $C_V$ is at most a $\Pn 5$ by (1.4).
If $L$ is a $\Pn 4$ and $T$ is a curve,
then the intersection $L\cap G$ is proper and $C_V$ has degree at most 5
and is rational or elliptic.\par If $T$ is a surface, then $T$ is a cubic
scroll
or a complete intersection of two quadrics.  In the latter case $T$ is a
cone or
a del Pezzo surface.\par
  If $T$ is neither a curve nor a surface in a $\Pn 4$, then it must be a
quadric hypersurface.\par
If $L$ is a $\Pn 3$, then $T$ is a curve of degree at most 4, or $T$ is a
quadric or $T$ is all of $L$.\par
If $L$ is contained in a plane, then $T$ is the whole plane, a conic or a
line. $\slutt$\medskip

This exhausts the list of possibilities for $T$.
Lemma (3.1) simplifies the analysis of each case.
Thus if $T$ is a quadric, then all the plane conics of $T$
correspond to quadrics of rank 4.  Now if two conics in $T$ meet in two
points,
then the corresponding quadrics of rank 4 have a common vertex.  By
choosing
different conics, we may conclude that the planes corresponding to the
points
on $T$ all have a common point, i.e. $V$ is a cone.\par
If $T$ is a (possibly degenerate) del Pezzo surface, we may again find a
conic
on it such that the hyperplanes through it cuts out a pencil of conics on
$T$.
The preceding argument shows that $V$ is a cone also in this case.\par
If $T$ is a cone over an elliptic quartic curve, then $C_V$ has degree
$\delta
=4\alpha +1$ or $\delta =4\alpha $ for some $\alpha$ depending on whether
$C_V$
meets the vertex of $T$ or not.  The corresponding genera are given by
$2g-2=2(2\alpha +1)(\alpha -1)$ and $2g-2=4\alpha (\alpha-1)$ respectively.
Combined with the inequality $\delta\leq 31$ from (2.2) it is easily
checked
that (1.1c) has a numerical solution only if $\alpha =1$, i.e. when $C_V$
is
elliptic. \par If $T$ is a $\Pn 3$, then the planes corresponding to the
points
of $T$ must all lie in a $\Pn 3$, so $V$ and $S$ is degenerate.\par
If $T$ is a quadric, then the argument above applies to show that $V$ is a
cone. \par
If $T$ is a curve in $\Pn 3$, then $C_V$ is rational or elliptic.
 If $T$ is contained in a plane, then $V$ is either contained in a $\Pn 3$
or it is a cone.  Combined with (1.5) we have shown that
\proclaim (3.3) Lemma.  $C_V$ is rational or elliptic, or it lies on a
cubic
scroll in a $\Pn 4$ or it lies on a quartic surface in $\Pn 5$ or $V$ is a
cone.\par \bigskip
{\bf  4  The cone case}\medskip
\proclaim (4.1) Proposition.  If $V$ is a cone, then it is a quadric or a
$\Pn
3$.\par
$Proof.$  If $V$ is a cone then there is some curve on $W$ which is
contracted
by the projection into $\Pn 4$.  Since it is contracted the numerical class
of this curve must be a multiple of the class $h^2-\delta h\cdot f$, where
$h$
is the calss of a hyperplane section and $f$ is the class of a fibre.
Unless
$C_V$ is rational (i.e. of degree $\delta \leq 2$), $S$ meets this curve in
at
most one point, so $0\leq S\cdot (h^2-\delta h\cdot f)\leq 1$.  Thus
$$0\leq
d-2\delta\leq 1.$$ If $d=2\delta$, then (1.1) yields
 $$2g-2={1\over 2}\delta^2-2\delta,$$
while if $d=2\delta +1$, then (1.1) yields
$$2g-2={1\over2}\delta^2-{3\over 2}\delta-1.$$
In both cases a comparison with Castelnuovos bound for the genus of space
curves
and with the genus of plane curves
shows that $\delta\leq 2$.
 $\slutt$\bigskip

{\bf 5 The cubic scroll case}\medskip
Assume that $C_V$ lies on a cubic scroll.  Let $E$ be a hyperplane section
of
the scroll and let $F$ be a member of the ruling, then numerically
$C_V\equiv \alpha E+ \beta F$  where $\alpha \geq 0$ and $\beta \geq
-\alpha$
 and $\delta =3\alpha +\beta$.  Since $C_V$ is smooth
we get by adjunction $2g-2=3{\alpha }^2-5\alpha +2\alpha\beta -2\beta$
 even if $T$ is singular.
With the inequality $\delta\leq 31$ of (2.2), a simple program checks the
possibilities and allow us to conclude that (1.1c)
 has no integral solutions in the given range.
Thus \proclaim (5.1) Proposition.  $C_V$ is not a cubic scroll.\par
\bigskip
{\bf 6 The quartic surface in $\Pn 5$ case.}\medskip
First assume that $C_V$ lies on a rational quartic scroll.
Let $E$ be a hyperplane section of the scroll
 and let $F$ be a member of the ruling, then numerically
$C_V\equiv \alpha E+ \beta F$
where $\alpha \geq 0$ and $\beta \geq -2\alpha$
 and $\delta =4\alpha +\beta$.  Since $C_V$ is smooth
we get by adjunction $2g-2=4{\alpha }^2-6\alpha +2\alpha\beta -2\beta$
 even if the scroll is singular.
With the inequality $\delta\leq 31$ of (2.2), a simple program checks the
possibilities and allow us to conclude that (1.1)
 has only two integral solutions in the given range.
Thus \proclaim (6.1) Lemma.  If $C_V$ lies on a rational quartic scroll,
then its numerical class is $3E-F$ or $6E+2F$.\par
If the quartic surface is not a scroll then it is a Veronese surface.
 So $C_V$ is a plane curve of degree $a$
embedded by conics in $\Pn 5$.  Thus $\delta =2a$ and $2g-2=a(a-3)$.
As above the relation (1.1c) is checked for $\delta \leq 31$ i.e. $a\leq
15$,
and there are no integral solutions. We have shown

\proclaim (6.2) Lemma. If $C_V$ spans $\Pn 5$, then it does not lie on a
Veronese surface.\par

We want to exclude the possibilities of (6.1) by a geometric argument:
First let $C_V$ be of numerical type  $3E-F$ on a scroll in the
grassmannian.
Then the degree of the conic bundle $S$ is 15 and the  sectional genus is
19.
Now, by (3.1), for any general line in the ruling of the scroll there is a
hyperplane section of  $S$ consisting of three conic sections and a
residual curve $A$.  The family of lines is rational so by Bertini the
curve
$A$ is irreducible for a general line in the ruling. The degree of $A$ is 9
while the arithmetic genus is 16. This is impossible. Similarly let  $C_V$
be of numerical type  $6E+2F$ on a scroll in the grassmannian.
 Then the degree of the conic bundle $S$ is 36 and the  sectional genus is
139.
For a general line in the ruling of the scroll there is a hyperplane
section of
$S$ consisting of six conic sections and an irreducible residual curve $A$.
The
degree of $A$ is 24 while the arithmetic genus is 133. This is impossible,
so we
may conclude \proclaim (6.3) Proposition. $C_V$ does not span $\Pn 5$.\par
\bigskip {\bf 7 Conclusion}\medskip

Combining (3.3), (4.1), (5.1) and (6.3) we are left with case that $C_V$ is
rational or elliptic. But by (1.4) this means that $2\leq\delta \leq 5$.
The
relations in (1.1) leave us with the possibility that $V$ is a quadric and
the
surface has degree 4 or 5 or that $V$ is a quartic and the surface is a
conic bundle of degree 8 over an elliptic curve.  The latter possibility
was
excluded by Okonek (cf. [Ok]).  Therefore we have
\proclaim (7.1) Theorem (Ellia,
Sacchiero).  The degree of a conic bundle in $\Pn 4$ is 4 or 5. \par
\proclaim
(7.2) Remark.  There are surfaces with a 1-dimensional family of conic
sections
which are not conic bundles, these surfaces are easily seen to be rational
and
are the cubic scrolls and the Veronese surfaces.

\bigskip
{\bf References}\bigskip

{\parindent 20pt
 \item {[A]} Aure, A.: On surfaces in projective 4-space.  Thesis, Oslo
(1987)
 \item{[BF]} Braun R., Fl\o ystad G.:  A bound for the degree of smooth
surfaces in $\Pn 4$ not of general type.   Preprint Bergen, Bayreuth (1993)
\item {[ES]} Ellia, P., Sacchiero, G. (to appear)
\item{[GP]} Gruson, L.,
Peskine, C.:  Genre des courbes de l'espace projectif. Algebraic Geometry
Proceedings, {Troms\o } 1977, SLN {\bf 687} (1978), 31-59.
 \item{[HJ]} Harris, J.:Curves In projective space, Le presses de
l'Universit\'e
de Montreal. (1982)
 \item{[HR]} Hartshorne, R.:  Algebraic Geometry.  Berlin,
Heidelberg, New York:  Springer (1977).
\item {[K]} Koelblen, L.: Surfaces de $\Pn 4$ trac\'ees sur une hyperface
cubique.  Journal f\"ur die Reine und Angewandte Mathematik, {\bf 433}
(1992),
113-141.
\item{[Ok]} Okonek, C.: Fl\"achen vom Grad 8 in $\Pn 4$. Math. Z. {\bf
191},
207-223 (1986).
\item{[Ra]} Ranestad, K.:  On smooth plane curve fibrations in $\Pn 4$.
Geometry of Complex Projective Varieties Proceedings. Cetraro 1990,
Mediterrenean Press, Rende, Italy, (1993)
 \item{[Ro]} Roth, L.: On the projective classification of
surfaces.  Proc. London Math. Soc. {\bf 42}  (1937) 142-170

\item{}}\bigskip
addresses:\par
Robert Braun\par
Mathematisches  Institut\par
Universit\"at Bayreuth\par
D-8580 Bayreuth\par
Germany\par
\medskip
Kristian Ranestad\par
Matematisk Institutt\par
boks 1053\par
N-0316 Oslo 3\par
Norway

\bye